\documentclass[10pt]{article}  
\usepackage{amsmath}
\usepackage{amssymb}
\usepackage{epsfig}
\usepackage{theorem}% (il faut utiliser le style theorem)

\theoremheaderfont{\scshape}
\newtheorem{thm}{Theorem}[section]

\theorembodyfont{\normalfont}

\numberwithin{equation}{section}

\title{Modeling the Coastal Ocean over a Time Period of Several Weeks}
\author{
P. Ailliot\thanks{Sabres, 
Universit\'e de Bretagne Sud, Centre Yves Coppens, Campus de Tohannic,
F-56000, Vannes} \and 
E. Fr\'enod\thanks{LMAM et Lemel, 
Universit\'e de Bretagne Sud, Centre Yves Coppens, Campus de Tohannic,
F-56000, Vannes}\and  
V. Monbet \thanks{Sabres et Lemel, 
Universit\'e de Bretagne Sud, Centre Yves Coppens, Campus de Tohannic,
F-56000, Vannes}
}

\newcommand{\nit}{\mathbb{N}}
\newcommand{\rit}{\mathbb{R}}

\newcommand{\eps}{{\varepsilon}}

\newcommand{\tc}{{t}}
\newcommand{\xc}{{x}}

\newcommand{\xvec}{{\mathbf{x}}}

\newcommand{\uvec}{\mathbf{u}}
\newcommand{\uuvec}{\mathbf{U}}

\newcommand{\eede}{{E}}
\newcommand{\hhde}{{\cal H}}
\newcommand{\hhdet}{\tilde \hhde}

\newcommand{\hhdemoy}{ \overline{\overline\hhde} }
\newcommand{\phhde}{{\cal I}}
\newcommand{\phhdet}{\tilde \phhde}

\newcommand{\vtot}{\mathbf{m}}
\newcommand{\vmr}{ \mathbf{M}}
\newcommand{\vmrt}{\tilde \vmr}

\newcommand{\vmrmoy}{\overline{\overline\vmr}}
\newcommand{\vmrunmoy}{ \overline{\overline{\vmr_1}} }
\newcommand{\vmrdemoy}{ \overline{\overline{\vmr_2}} }
\newcommand{\pvmr}{\mathbf{N}}
\newcommand{\pvmrt}{\tilde \pvmr}

\newcommand{\vwn}{\mathbf{W}}
\newcommand{\vwnt}{\tilde \vwn}

\newcommand{\vwnmoy}{\overline{\overline{\vwn}}}
\newcommand{\vwnunmoy}{\overline{\overline{\vwn_1}}}
\newcommand{\vwndemoy}{\overline{\overline{\vwn_2}}}
\newcommand{\fxt}{\mathbf{F}}

\newcommand{\tov}{{\overline{t}}}
\newcommand{\llov}{{\overline{L}}}
\newcommand{\lov}{{\overline{l}}}

\newcommand{\eedeov}{{\overline{E}}}
\newcommand{\hhdeov}{{\overline{H}}}
\newcommand{\phhdeov}{{\overline{I}}}
\newcommand{\vmrov}{{\overline{M}}}
\newcommand{\pvmrov}{{\overline{N}}}
\newcommand{\vwnov}{{\overline{W}}}
\newcommand{\fxtov}{{\overline{F}}}
\newcommand{\omov}{{\overline{\omega}}}

\newcommand{\ds}{\displaystyle}

\newcommand{\tsep}{\frac{t}{\eps}}
\newcommand{\unsep}{\frac{1}{\eps}}

\newcommand{\fracp}[2]{\frac{\partial #1}{\partial #2}}
\newcommand{\fracpse}[2]{\frac{\partial ^2 #1}{\partial #2 ^2}}
\newcommand{\fracpsed}[3]{\frac{\partial ^2 #1}{\partial #2 \partial #3}}

\begin{document}
\maketitle

{\small {\bf Abstract - } From a scale analysis of hydrodynamic phenomena
having a significant action on the drift of an object in coastal ocean waters, 
we deduce equations
modeling the associated hydrodynamic fields over a time period of several weeks.
These models are essentially non linear hyperbolic systems of PDE involving 
a small parameter.
Then from the models we extract a simplified and nevertheless typical one
for which we prove that its classical solution exists on a time interval
which is independent of the small parameter. We then show that the solution weak$-*$ converges as the small parameter goes to zero and we characterize the equation satisfied by the weak$-*$ limit.
}

{\small {\bf Keywords - } Hyperbolic PDE, coastal ocean modeling, scale analysis, asymptotic analysis
}

%\tableofcontents

\section{Introduction}
This paper is part of a work program concerning the modeling of object drift
in near coastal ocean waters over a several week time period.

The final target of this program is to develop methods to forecast the drift of things
like containers, lost objects or oil spill over long periods of time in near coastal
ocean areas. Such methods would be of interest for services in charge of maritime
safety, environmental studies or pollution impact assessment. To reach this target,
several research topics need to be further investigated. For instance, improvements are needed
in the field of the numerical methods to simulate long term drift, in the modeling
and simulation of the near coastal ocean waters, in the understanding of ocean-object and
ocean-spill interactions and, of course, in the integration of all those aspects 
to move toward a coherent theory.

~

In the previous paper of this work program, 
Ailliot, Fr\'enod and Monbet \cite{AilFreMon} considered the numerical method facet.
We built a numerical method coupling a two scale expansion method,
explored in  \cite{frenod2004},
and a stochastic wind simulator,
in the spirit of \cite{AilliotThese,Ailliot2006,MonbetAlliotPrevosto},
in order to estimate probability of events that may happen to the considered object,
such as running aground in a given area.
In \cite{AilFreMon}, the simplified model was supposed to describe
the object dynamic in the ocean. It involved an ocean velocity
field which was decomposed into a sum of a velocity due to the tide wave and of a 
perturbation. Both of them were periodic of the tide period and with modulated amplitude and
the fields used for the numerical simulations were not realistic.

~

The present paper deals with the modeling of near coastal ocean.
The sea velocity and the fluctuation of the sea level due to the tide wave are 
well known in many coastal areas of the world.
%: Data are placed at the disposal by state services in charge of maritime safety which allow, by simple interpolation rules, a precise computation of those fields at any place and any time.
Those main fields are perturbed by fields with a smaller order of magnitude having 
a net long-term result. Such perturbations, which are produced by meteorological 
factors, propagate and interact with the main fields.
The precise aim of the present paper is to make a first step toward the
set up of a modeling procedure in order to establish
partial differential equation systems describing the evolution 
of those perturbations and to suggest ideas to solve them.  

~

The paper is organized as follows. In section \ref{Res}, we summarize the main mathematical results.
 
Then, in section \ref{ModSec}, we set up the previously evoked modeling procedure.
The models which are built in this section are deduced from the Shallow Water Equations via a scale analysis
of the geophysical phenomena concerned and of the geometrical size of the concerned 
domain.
As for the size of the coastal domain and the characteristic order of magnitude
of the wind velocity, we consider several possibilities 
giving rise to several models. They are all essentially hyperbolic systems of 
partial differential equations
with a singular perturbation involving a small parameter.

From those models, in section \ref{ExisSec}, we extract a simplified 
one which is nevertheless
typical. For it, adapting classical methods for hyperbolic systems (see
Kato \cite{Kato:1975}, Majda \cite{Majda:1984},
Klainerman and Majda \cite{KlaiMaj:1981,KlaiMaj:1982},
Schochet \cite{schochet:1986, schochet:1988},
and Metivier and Schochet \cite{MetScho}), 
we prove that its classical solution exits with a time existence
independent of the small parameter. 

In section \ref{AsBe}, using a homogenisation method
(see Tartar \cite{tartar:1977},
Bensoussan, Lions and Papanicolaou \cite{bensoussan/lions/papanicolaou:1978},
Sanchez-Palencia \cite{sanchez-palencia:1978} and 
Lions \cite{lions:PS}),
we set that this classical solution  weak$-*$ converges to a function.
Keeping within the mind frame of Fr\'enod \cite{frenod:1994},
Fr\'enod and Hamdache \cite{frenod/hamdache:1996}, 
Fr\'enod, Raviart and Sonnendr\"ucker \cite{FRS:1999},
Joly, M\'etivier and Rauch \cite{JMR93} or 
Schochet \cite{schochet:1994}
we finally look for the form of this function
and establish the equations allowing for its computation.

Finally, in section \ref{conclu} we conclude and give some perspectives.
 
~

{\bf acknowledgments -} The authors want to thank Joanna Roppers for proofreading
the manuscript.

\section{Results}
\label{Res}
In this section we present the main results. We first present one of the
models, involving a small parameter, set out in this paper. Then we state
a Theorem claiming the existence of the solution to a simplified
version of this model. Finally, we exhibit the asymptotic behavior 
of this solution as the small parameter goes to zero.  

~

The model we present now, and which is set out among others in section \ref{ModSec},
describes the evolution, over a several months time period, of the perturbation of the sea velocity and of the sea
level in an ocean domain above the continental shelf at a latitude about $45^o$
and with stormy weather conditions.

The small parameter involved in this model is the ratio tide duration 
on observation time scale.
The first one is about 13 hours and the second is about three months. Hence 
the involved small parameter is:
\begin{equation}\label{PS2.0}
\eps = \frac{1}{200}.
\end{equation}
Variables and fields involved are all rescaled; rescaled meaning that the order 
of magnitude of those variables and fields is one and that they have no 
physical dimension.

The rescaled velocity of the sea $\vmrt$ and the water depth $\hhdet$ 
induced by the tide wave are considered as known
and periodic with modulated amplitude. In other words, $t$ being the rescaled time and
$\xvec$ the rescaled position,
\begin{equation}
\label{hypdeptheta}
\vmrt(t,\xvec)= \vmr(t,\tsep,\xvec) \text{ and }
\hhdet(t,\xvec)= \hhde(t,\tsep,\xvec),
\end{equation}
where  $\vmr$ and $\hhde$ are regular functions and where 
$\theta\mapsto (\vmr(t,\theta,\xvec),\hhde(t,\theta,\xvec))$ 
is $1-$periodic.

The model says that 
the total sea velocity, expressed in $km/h$, writes $0.5(\vmrt+\eps\pvmrt)$ 
and that the total sea level is
$\frac{3}{2\eps}(\eede+2\eps\hhdet+2\eps^2\phhdet)$, where $\eede$ is the rescaled mean 
sea level and where $\pvmrt$ and $\phhdet$ are rescaled perturbations.
Moreover, $(\phhdet, \pvmrt)$ is solution to
\begin{align} 
\ds
&
\begin{aligned}
  \fracp{\phhdet}{\tc} + 
  \nabla (\frac{1}{\eps}\eede+2\hhdet)\cdot\pvmrt+
  (\frac{1}{\eps}\eede+2\hhdet)\nabla\cdot\pvmrt+
  2(\nabla \phhdet)\cdot\vmrt +
  2\phhdet (\nabla\cdot\vmrt) &
\\ \ds + 
  2\eps\big( (\nabla \phhdet)\cdot\pvmrt+
  \phhdet(\nabla\cdot\pvmrt)\big)&=0, 
\end{aligned}
\label{REcs3intr}
\\ \ds
&
\begin{aligned}
\fracp{\pvmrt}{\tc} +
  2(\nabla\pvmrt)\vmrt+
  2(\nabla\vmrt)\pvmrt+
  2\eps (\nabla\pvmrt)\pvmrt+
  \frac{\pi}{2\eps}\pvmrt^\perp+
  \frac{1}{4\eps}\nabla\phhdet-
  13\eps^4\Delta  \vmrt-
  13\eps^5 \Delta  \pvmrt ~~~~&
\\
  -13\eps^4
  \frac{(\nabla \vmrt)\nabla(\eede+2\eps\hhdet)}
       {\eede+2\eps\hhdet +2\eps^2\phhdet}
 -26\eps^6
     \frac{(\nabla\vmrt)\nabla\phhdet }
          {\eede+2\eps\hhdet+2\eps^2\phhdet}
 -13\eps^5
   \frac{(\nabla\pvmrt ) \nabla(\eede+2\eps\hhdet)}
        {\eede+2\eps\hhdet+2\eps^2\phhdet}&
\\
 -26\eps^7
   \frac{(\nabla\pvmrt ) \nabla\phhdet}
        {\eede+2\eps\hhdet+2\eps^2\phhdet}&
\\
+\frac3\eps \frac{\frac{1}{\eede+2\eps\hhdet+ 2\eps^2\phhdet}}
     {1+\frac{0.8}{\eps^2}(\eede+2\eps\hhdet+2\eps^2\phhdet)} \vmrt
+3\frac{\frac{1}{\eede+2\eps\hhdet+ 2\eps^2\phhdet}}
     {1+\frac{0.8}{\eps^2}(\eede+2\eps\hhdet+2\eps^2\phhdet)}  \pvmrt
 &= 
\\ 
   6\frac{\frac{1}{\eede+2\eps\hhdet+2\eps^2\phhdet}}
        {1+\frac{1.5}{\eps}
        (\eede+2\eps\hhdet+2\eps^2\phhdet)} 
        (\frac1\eps\vwnt- \vmrt)- 
   6\eps\frac{\frac{1}{\eede+2\eps\hhdet+2\eps^2\phhdet}}
        {1+\frac{1.5}{\eps}
        (\eede+2\eps\hhdet+2\eps^2\phhdet)} &\pvmrt.
\end{aligned}
\label{REcs4intr}
\end{align}
In this system, $\vwnt$ is the rescaled wind velocity, 
$\pvmrt^\perp=(-\pvmrt_2,\pvmrt_1)$, 
$\Delta$ stands for the laplacian,
$\nabla\cdot$ for the divergence operator and 
$\nabla$ stands for the gradient of
scalar fields and for the Jacobian matrix of bi-dimensional fields.

~

Motivated by this system, we consider a simplified version of it 
which consists in considering that the ocean bottom is flat, i.e. 
$\eede\equiv 1$, in forgetting all the power of $\eps$ greater than $1$
and in setting all constants to 1:
\begin{align}  \ds
&
\begin{aligned}
  \fracp{\phhdet}{\tc} + 
  (\nabla \hhdet)\cdot\pvmrt+
  (\frac{1}{\eps}+\hhdet)(\nabla\cdot\pvmrt)+
  (\nabla \phhdet)\cdot\vmrt +
  \phhdet (\nabla\cdot\vmrt) + 
  \eps\big( (\nabla \phhdet)\cdot\pvmrt+
  \phhdet(\nabla\cdot\pvmrt)\big)&=0, 
\end{aligned}
\label{SM1cs3intr}
\\ \ds
&
\fracp{\pvmrt}{\tc} +
  (\nabla\pvmrt)\vmrt+
  (\nabla\vmrt)\pvmrt+
  \eps (\nabla\pvmrt)\pvmrt+
  \frac{1}{\eps}\pvmrt^\perp+
  \frac{1}{\eps}\nabla\phhdet
 =\vwnt.
\label{SM1cs2intr}
\end{align}
In this system, $t\in[0,T]$, $\xvec=(x_1,x_2)\in\rit^2$. The unknowns are 
$\phhdet\equiv \phhdet(t,\xvec)$ and $\pvmrt\equiv \pvmrt(t,\xvec)$.
Their evolution is influenced by  $\vmrt$ and $\hhdet$ for which we 
assume (\ref{hypdeptheta}) and by $\vwnt$ for which we also assume
\begin{equation}
\label{hypdepthetawnd}
\vwnt(t,\xvec)= \vwn(t,\tsep,\xvec) ,
\end{equation}
with function $\vwn$ regular and with $\theta\mapsto \vwn(t,\theta,\xvec)$ 
being  $1-$periodic.
This assumption is not really convenient for real wind 
velocity time series but is comfortable from a mathematical point of view (see Ailliot, Fr\'enod and Monbet \cite{AilFreMon} for a more detailed discussion).
Moreover, we equip this system with the following initial conditions
\begin{gather}
\label{Cini}
\phhdet_{|t=0} = \phhdet_0, ~~~~\pvmrt_{|t=0} = \pvmrt_0,
\end{gather}
and we can claim the following Theorem.
\begin{thm}
\label{THM1}
Under assumptions (\ref{hypdeptheta}) and (\ref{hypdepthetawnd}),
if $(\phhdet_0,\pvmrt_0)\in (H^s(\rit^2))^3$ with $s>3$, then there exists
a time $T$, not depending on $\eps$, such that the classical solution 
$(\phhdet,\pvmrt)$ 
$\in (C([0,T],(H^s(\rit^2))^3 )\cap(C^1([0,T],(H^{s-1}(\rit^2))^3))$ 
of (\ref{SM1cs3intr}), (\ref{SM1cs2intr}) and (\ref{Cini}) exits and is unique.
Moreover this solution satisfies
\begin{equation}
\label{estimsol}
\sup_{t\in[0,T]}\|(\phhdet,\pvmrt)\|_s \leq c,
\end{equation}
for a constant c not depending on $\eps$, 
where $\|~\|_s$ stands for the norm in $(H^s(\rit^2))^3$.
\end{thm}

Concerning the asymptotic behavior of $(\phhdet,\pvmrt)$, as $\eps$ goes to zero,
 we have the following result.
\begin{thm}
\label{THM2}
Under the assumptions of Theorem \ref{THM1}, there exist functions
$\phhde\equiv \phhde(t,\xvec)$ 
$\in C([0,T],$ $H^s(\rit^2))$
%\cap(C^1([0,T]\times\rit,H^{s-1}(\rit^2)))
and $\pvmr\equiv \pvmr(t,\xvec)$
$\in C([0,T],(H^s(\rit^2))^2)$,
%\cap(C^1([0,T]\times\rit,H^{s-1}(\rit^2)))
such that as $\eps$ goes to $0$, the solution $(\phhdet,\pvmrt)$
of (\ref{SM1cs3intr}), (\ref{SM1cs2intr}) and (\ref{Cini}) 
weak$-*$ converges to $(\phhde ,\pvmr )$ in $L^\infty([0,T],$ $(H^s(\rit^2))^3)$.
Moreover, $\phhde$ and $\pvmr$  are linked by 
\begin{equation}
\label{ContrCons}
%\phhde(t,\xvec)=  g(t,\xvec), ~~~
\pvmr_1(t,\xvec)= - \fracp{\phhde}{x_2} (t,\xvec) , ~~~
\\\pvmr_2(t,\xvec)= \fracp{\phhde}{x_1} (t,\xvec) ,
\end{equation}
and  $\phhde$ is solution to 
\begin{multline}
\label{SysEff}
\fracp{\ds \big( \phhde-\Delta \phhde\big)}{t} + \vmrmoy  \cdot \nabla  \phhde
-\fracp{\ds \Big(\vmrunmoy \fracpse{\phhde}{x_1}\Big)}{x_1}
-\fracp{\ds \Big(\vmrdemoy \fracpsed{\phhde}{x_1}{x_2}\Big)}{x_1}
\\
-\fracp{\ds \Big(\vmrunmoy \fracpsed{\phhde}{x_1}{x_2}\Big)}{x_2}
-\fracp{\ds \Big(\vmrdemoy \fracpse{\phhde}{x_2}\Big)}{x_2} 
-(\nabla\hhdemoy)^\perp \cdot \nabla \phhde + (\nabla\cdot\vmrmoy) \phhde
\\
+ \fracp{\ds \Big(\fracp{\vmrdemoy}{x_1} \fracp{\phhde}{x_2} \Big)}{x_1}
- \fracp{\ds \Big(\fracp{\vmrdemoy}{x_2} \fracp{\phhde}{x_1} \Big)}{x_1}
- \fracp{\ds \Big(\fracp{\vmrunmoy}{x_1} \fracp{\phhde}{x_2} \Big)}{x_2}
+ \fracp{\ds \Big(\fracp{\vmrunmoy}{x_2} \fracp{\phhde}{x_1} \Big)}{x_2}
\\
= \fracp{\vwnunmoy}{x_2}  - \fracp{\vwndemoy}{x_1},
\end{multline}
where $\ds \vmrmoy =\int_0^1 \vmr \, d\theta$,
$\ds \hhdemoy =\int_0^1 \hhde \, d\theta$ and
$\ds \vwnmoy=\int_0^1 \vwn \, d\theta$, and 
equipped with initial conditions
\begin{gather}
\label{SysEffCI}
\big( \phhde-\Delta \phhde\big)_{|t=0} 
= \phhdet_0 + \fracp{(\pvmrt_0)_1}{x_1} - \fracp{(\pvmrt_0)_2}{x_2}.
\end{gather} 
\end{thm}
\section{Models}
\label{ModSec}
In this section, we first consider a reference model. It
consists in removing the ocean level and the 
ocean velocity which are induced by the tide wave
from the Shallow Water Equations. This gives rise
to a system of equations governing the time evolution of the ocean level perturbation
and of the ocean velocity perturbation. Then, we analyse the scale of the variables
and fields involved in the problem we want to describe. Rescaling the reference
model in view of this scale analysis finally yields the desired models.

\subsection{Reference model}
It is generally admitted that 
the evolution of the ocean level $h\equiv h(\tc,\xvec)$
and of the ocean velocity $\vtot \equiv\vtot(\tc,\xvec)$
is well described by the following Shallow Water Equations
\begin{align} \ds
&\fracp{h}{\tc} + \nabla (h-h_b) \cdot \vtot + (h-h_b) \nabla\cdot \vtot =0,
\label{StV1}
\\  \ds
&\begin{aligned}
\fracp{\vtot}{\tc} +\big(\nabla\vtot\big)\vtot
+ f \vtot^\perp 
+ g \nabla h - c \Delta \vtot -c\frac{\big(\nabla\vtot\big)\nabla(h-h_b)}{h-h_b}
+ \frac{\frac{\kappa}{h-h_b}}{1+\frac{\kappa}{c}(h-h_b)}\vtot &=
\\  \ds
\frac{\frac{\mu}{h-h_b}}{1+\frac{\mu}{c}(h-h_b)} (\vwnt -\vtot) &+\fxt,
\end{aligned}
\label{StV2}
\end{align}
equipped with ad-hoc initial and boundary conditions.
This system was introduced by 
Saint-Venant \cite{saintvenant:1871}.
For an exhaustive explanation concerning ocean modeling 
and the construction of this model we refer for instance to
Pedlosky \cite{pedlosky},
Nihoul \cite{Nihoul},
Lions, Temam and Wang \cite{LionsTemamWang1},
Stoker \cite{stoker},
Whitham \cite{whitham} or
Johnson \cite{johnson}.
For a deduction of the Shallow Water Model taking into account
viscosity, being able to model the consequences of wind and bottom 
friction actions, which is considered here, we refer to  
Gerbeau and Perthame \cite{GerPer}.
In system (\ref{StV1})-(\ref{StV2}), $h_b\equiv h_b(\xvec)$ is the depth of 
the ocean bottom, $f$ is the Coriolis parameter, $g$ is the gravity acceleration
and $c$ is the water viscosity. 
The friction coefficient on the bottom is $\kappa$ and the 
air-water friction coefficient is $\mu$. 
Lastly, $\vwnt\equiv\vwnt(\tc,\xvec)$ is the wind velocity and 
$\fxt$ may take into account the action of other meteorological factors
like atmospheric pressure.

~

Now we isolate the action of the tide wave.
In other words, we consider that the ocean depth variation 
$\hhdet\equiv\hhdet(t,\xvec)$ around the mean water height
$\eede\equiv\eede(\xvec)$ and the ocean velocity 
$\vmrt\equiv\vmrt(t,\xvec)$ which are induced by the tide wave
are known. We consider that $(\eede+\hhdet,\vmrt)$ is the solution to
\begin{align} \ds
&\fracp{\hhdet}{\tc} + \nabla (\eede+\hhdet) \cdot \vmrt + 
(\eede+\hhdet) \nabla\cdot \vmrt =0,
\label{EqTW1}
\\  \ds
&\fracp{\vmrt}{\tc} +\big(\nabla\vmrt\big)\vmrt
+ f \vmrt^\perp 
+ g \nabla(\eede+\hhdet+h_b)  =0,
\label{EqTW2}
\end{align}
\begin{figure}
\centering \epsfig{file=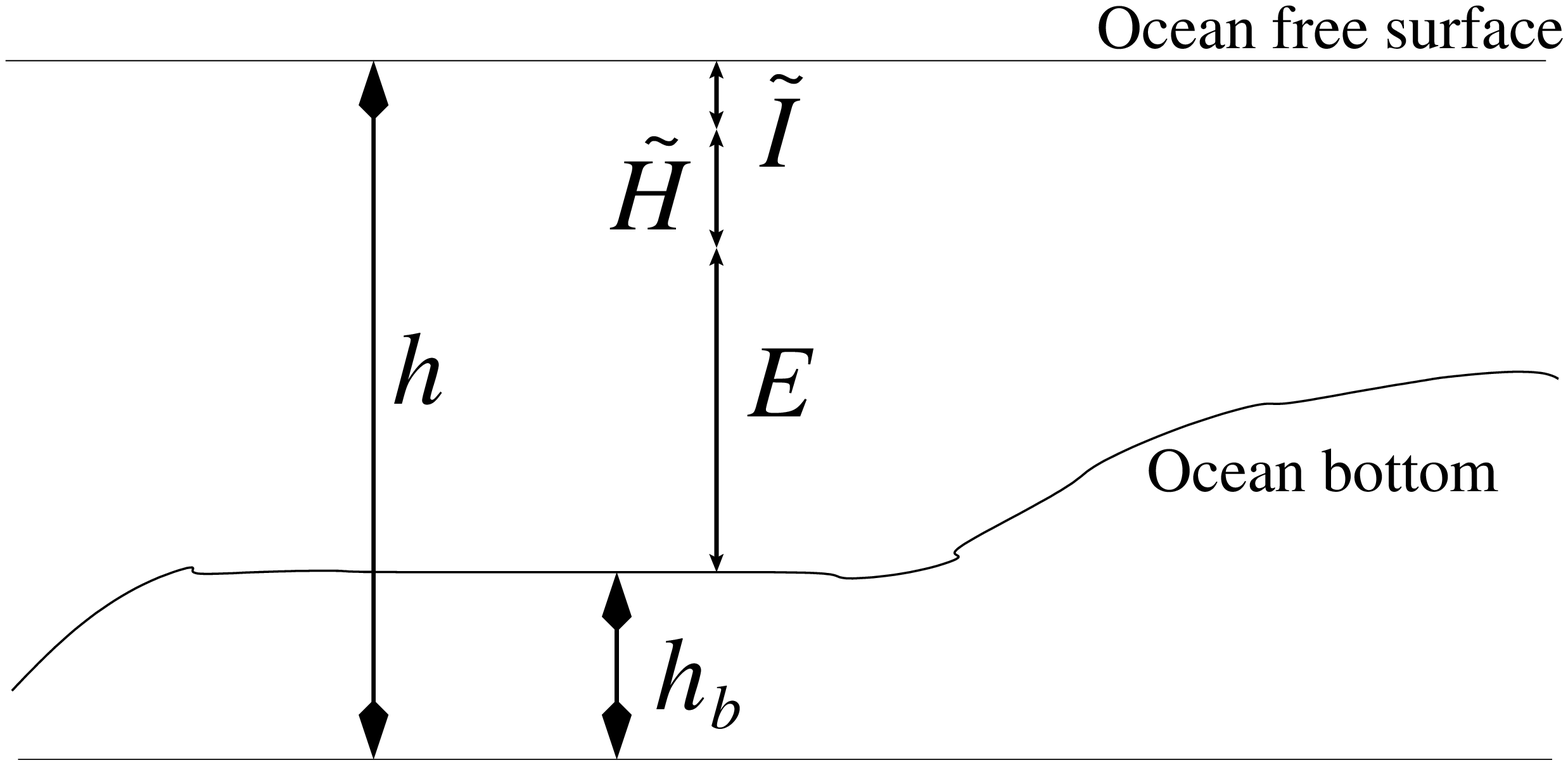,height=4cm}
\caption{Fields $h, h_b,\eede, \hhdet$ and $\phhdet $}\label{fig1}
\end{figure}
with initial and  boundary conditions imposed by the 
tide wave.

A brief parameter size analysis, that will be confirmed in the 
next sections, shows that the terms 
$- c \Delta \vtot -c\frac{\big(\nabla\vtot\big)\nabla(h-h_b)}{h-h_b}
+ \frac{\frac{\kappa}{h-h_b}}{1+\frac{\kappa}{c}(h-h_b)}\vtot$
have a very small influence on the sea movement. Hence we have chosen
to put them in the equations for the perturbations hereafter. 

~

Now we introduce the perturbations $\phhdet$ and $\pvmrt$ which are
defined such that $h=h_b+\eede+\hhdet+\phhdet$ and $\vtot=\vmrt+\pvmrt$   
(see Figure \ref{fig1}).

Replacing $h$ and $\vtot$ by these expressions in (\ref{StV1})-(\ref{StV2})
and removing the terms appearing in (\ref{EqTW1})-(\ref{EqTW2}) leads to the equations for
$\phhdet$ and $\pvmrt$. Then we obtain the following reference system
which is the starting point of our scale analysis.
\begin{align} \ds
&\fracp{\phhdet}{t} + 
\nabla (\eede+\hhdet))\cdot\pvmrt+
(\eede+\hhdet)(\nabla\cdot\pvmrt)+
(\nabla \phhdet)\cdot\vmrt +
\phhdet (\nabla\cdot\vmrt) + 
(\nabla \phhdet)\cdot\pvmrt+
\phhdet(\nabla\cdot\pvmrt)=0, 
\label{mtl3}\\ \ds
&\begin{aligned}
\fracp{\pvmrt}{t} +
(\nabla\pvmrt)\vmrt+
(\nabla\vmrt)\pvmrt+
(\nabla\pvmrt)\pvmrt
+f\pvmrt^\perp +
g \nabla\phhdet - 
c \Delta \vmrt-
c \Delta  \pvmrt ~~~~~~~~~~~~~~~~~~& 
\\
-c\frac{\big(\nabla\vmrt\big)\nabla(\eede+\hhdet) }{\eede+\hhdet+\phhdet}
-c\frac{(\nabla\vmrt)\nabla\phhdet }{\eede+\hhdet+\phhdet}
-c\frac{(\nabla\pvmrt)\nabla (\eede+\hhdet)}{\eede+\hhdet+\phhdet}
-c\frac{(\nabla\pvmrt)\nabla \phhdet}{\eede+\hhdet+\phhdet}&
\\
+\frac{\frac{\kappa}{\eede+\hhdet+\phhdet}}
     {1+\frac{\kappa}{c}(\eede+\hhdet+\phhdet)} \vmrt
+\frac{\frac{\kappa}{\eede+\hhdet+\phhdet}}
     {1+\frac{\kappa}{c}(\eede+\hhdet+\phhdet)}  \pvmrt &=
\\
\frac{\frac{\mu}{\eede+\hhdet+\phhdet}}
     {1+\frac{\mu}{c}(\eede+\hhdet+\phhdet)} (\vwnt - \vmrt)
-\frac{\frac{\mu}{\eede+\hhdet+\phhdet}}
     {1+\frac{\mu}{c}(\eede+\hhdet+\phhdet)} & \pvmrt + \fxt.
\label{mtl4}
\end{aligned}
\end{align}

\subsection{Rescaled variables and fields}
We introduce a reference time $\tov$, two reference lengths  
$\llov$ and  $\lov$. Those reference values, as well as the other ones introduced hereafter,
will represent characteristic values (mean or maximum values for example) of the
physical quantities under consideration.
We consider the rescaled variables $t'$ and $\xvec'=$($\xvec_1'$, 
$\xvec_2')$ expressing time and  position
in unit $\tov$, $\llov$ and   $\lov$.
They are defined as
\begin{equation}
t=\tov t', \xvec_1 =\llov\xvec_1' \text{ and  } \xvec_2 =\lov\xvec_2'.
\end{equation}
If the reference values are chosen as evoked above, the order of 
magnitude of the rescaled variables are 1.
Then we define $\vmrov$ and $\pvmrov$ as the characteristic velocity
of the tide wave and its perturbation; $\eedeov$ the characteristic
value of the mean water depth,
$\hhdeov$ the  characteristic 
tidal range and $\phhdeov$ the characteristic value of its perturbation.
$\vwnov$ is the characteristic wind velocity and  $\fxtov$
the  characteristic scale of the field $\fxt$.
The rescaled fields have the following definitions:
\begin{alignat}{2}\ds &
\vmrt'(t',\xvec') = \frac{1}{\vmrov} \vmrt (\tov t',\llov\xvec'_1,\lov\xvec'_2),
&~~~~~~
&
\pvmrt' (t',\xvec') = \frac{1}{\pvmrov} 
\pvmrt (\tov t',\llov\xvec'_1,\lov\xvec'_2),
\\ \ds & 
\eede'(\xvec') = \frac{1}{\eedeov}
\eede(\llov\xvec'_1,\lov\xvec'_2),&~~~~~~
&~~
\\ \ds &
\hhdet'(t',\xvec') = \frac{1}{\hhdeov}
\hhdet(\tov t',\llov\xvec'_1,\lov\xvec'_2),&~~~~~~
&
\phhdet'(t',\xvec') = \frac{1}{\phhdeov}
\phhdet(\tov t',\llov\xvec'_1,\lov\xvec'_2),~~~~~~
\\ \ds & 
\tilde\vwn'(t',\xvec') = \frac{1}{\vwnov} \vwn(\tov t',\llov\xvec'_1,\lov\xvec'_2),
&~~~~~~
&
\fxt'(t',\xvec') = \frac{1}{\fxtov} \fxt(\tov t',\llov\xvec'_1,\lov\xvec'_2).
\end{alignat}
Lastly, we introduce $\omov$ the tide wave frequency and we assume
that $\vmrt$ and  $\hhdet$ are  $1/\omov-$periodic function with modulated amplitude.
In other words, we assume
\begin{equation}
\vmrt'(t',\xvec')=  \vmr'( t',\omov\tov t',\xvec'), ~~~~~~~~
\hhdet'(t',\xvec')=  \hhde'( t',\omov\tov t',\xvec'),
\end{equation}
where  
$\theta\mapsto (\vmr'(t',\theta,\xvec'),\hhde'(t',\theta,\xvec'))$ 
is $1-$periodic.

From system (\ref{mtl3})-(\ref{mtl4}), we deduce the following rescaled equations
for $\pvmrt'$, and $\phhde'$ with known fields $\vmr'$, $\eede'$, $\vwn'$  and $\fxt'$:
\begin{align}\ds
&
\begin{aligned}
  \fracp{\phhdet'}{\tc'} &+ \frac{\hhdeov}{\phhdeov} \frac{\pvmrov\tov}{\llov}
  \Bigg[
  \begin{pmatrix}\ds\frac{\eedeov}{\hhdeov}\fracp{\eede'}{\xc'_1}+
    \fracp{\hhdet'}{\xc'_1}\\ \ds
    \frac{\llov}{\lov}(\frac{\eedeov}{\hhdeov}\fracp{\eede'}{\xc'_2}+
    \fracp{\hhdet'}{\xc'_2})
  \end{pmatrix}
  \cdot\pvmrt'+
  (\frac{\eedeov}{\hhdeov}\eede'+\hhdet')
  \big(\fracp{\pvmrt_1'}{\xc'_1}+\frac{\llov}{\lov}\fracp{\pvmrt_2'}{\xc'_2}\big)\Bigg] + 
  \frac{\vmrov\tov}{\llov}\Bigg[
  \begin{pmatrix}\ds\fracp{\phhdet'}{\xc'_1}\\ 
         \ds\frac{\llov}{\lov}\fracp{\phhdet'}{\xc'_2}
  \end{pmatrix}
  \cdot\vmrt' \\
  &+\phhdet'\big(\fracp{\vmrt'_1}{\xc'_1}+\frac{\llov}{\lov}\fracp{\vmrt'_2}{\xc'_2}\big)\Bigg]
  + \frac{\pvmrov\tov}{\llov}\Bigg[\begin{pmatrix}\ds\fracp{\phhdet'}{\xc'_1}\\
                             \ds\frac{\llov}{\lov}\fracp{\phhdet'}{\xc'_2}\end{pmatrix}
  \cdot\pvmrt'+
  \phhdet'\big(\fracp{\pvmrt'_1}{\xc'_1}+\frac{\llov}{\lov}\fracp{\pvmrt'_2}{\xc'_2}\big)\Bigg]  =0, 
\end{aligned}
\label{ppmtl3}
\end{align}
\begin{align}
&\begin{aligned}
  \fracp{\pvmrt'}{\tc'} +\frac{\vmrov\tov}{\llov}\Bigg[
  \Big( \fracp{\pvmrt'}{\xc'_1},\frac{\llov}{\lov}\fracp{\pvmrt'}{\xc'_2}\Big)\vmrt'+
  \Big( \fracp{\vmrt'}{\xc'_1},\frac{\llov}{\lov}\fracp{\vmrt'}{\xc'_2}\Big)\pvmrt'\Bigg]+
  \frac{\pvmrov\tov}{\llov}
       \Big( \fracp{\pvmrt'}{\xc'_1},\frac{\llov}{\lov}\fracp{\pvmrt'}{\xc'_2}\Big)\pvmrt'
  +f \tov \pvmrt'^\perp~~~~~~~~~~~~~&
\\
  +\frac{g\tov}{\pvmrov}\frac{\phhdeov}{\llov}\begin{pmatrix}
     \ds\fracp{\phhdet'}{\xc'_1}\\\ds\frac{\llov}{\lov}\fracp{\phhdet'}{\xc'_2}
   \end{pmatrix}
  - \frac{c\tov}{\llov^2}\frac{\vmrov}{\pvmrov}
     (\fracpse{\vmrt'}{{\xc'_1}}+\frac{\llov^2}{\lov^2}\fracpse{\vmrt'}{{\xc'_2}})
  - \frac{c\tov}{\llov^2}
     (\fracpse{\pvmrt'}{{\xc'_1}}+\frac{\llov^2}{\lov^2}\fracpse{\pvmrt'}{{\xc'_2}})&
\\
   -\frac{c\tov}{\llov^2}\frac{\vmrov}{\pvmrov}
  \frac{ \ds
  \Big( \fracp{\vmrt'}{\xc'_1},\frac{\llov}{\lov}\fracp{\vmrt'}{\xc'_2}\Big)
  \begin{pmatrix}\ds\fracp{\eede'}{\xc'_1}+
    \frac{\hhdeov}{\eedeov}\fracp{\hhdet'}{\xc'_1}\\ \ds
    \frac{\llov}{\lov}(\fracp{\eede'}{\xc'_2}+
    \frac{\hhdeov}{\eedeov}\fracp{\hhdet'}{\xc'_2})
  \end{pmatrix}
       }
       {\eede'+\frac{\hhdeov}{\eedeov}\hhdet' +\frac{\phhdeov}{\eedeov}\phhdet'}&
\\ 
   -\frac{c\tov}{\llov^2}\frac{\vmrov}{\pvmrov}\frac{\phhdeov}{\eedeov}
     \frac{\ds 
    \Big( \fracp{\vmrt'}{\xc'_1},\frac{\llov}{\lov}\fracp{\vmrt'}{\xc'_2}\Big)
    \begin{pmatrix}
     \ds\fracp{\phhdet'}{\xc'_1}\\\ds\frac{\llov}{\lov}\fracp{\phhdet'}{\xc'_2}
   \end{pmatrix}
          }
          {\eede'+\frac{\hhdeov}{\eedeov}\hhdet'
                                              +\frac{\phhdeov}{\eedeov}\phhdet'}
   -\frac{c\tov}{\llov^2}
   \frac{ \ds
     \Big( \fracp{\pvmrt'}{\xc'_1},\frac{\llov}{\lov}\fracp{\pvmrt'}{\xc'_2}\Big)
     \begin{pmatrix}\ds\fracp{\eede'}{\xc'_1}+
      \frac{\hhdeov}{\eedeov}\fracp{\hhdet'}{\xc'_1}\\ \ds
      \frac{\llov}{\lov}(\fracp{\eede'}{\xc'_2}+
      \frac{\hhdeov}{\eedeov}\fracp{\hhdet'}{\xc'_2})
     \end{pmatrix}
        }
               {\eede'+\frac{\hhdeov}{\eedeov}\hhdet'+\frac{\phhdeov}{\eedeov}\phhdet'}&
\\
   -\frac{c\tov}{\llov^2}\frac{\phhdeov}{\eedeov} 
   \frac{ \ds
     \Big( \fracp{\pvmrt'}{\xc'_1},\frac{\llov}{\lov}\fracp{\pvmrt'}{\xc'_2}\Big)
     \begin{pmatrix}
       \ds\fracp{\phhdet'}{\xc'_1}\\\ds\frac{\llov}{\lov}\fracp{\phhdet'}{\xc'_2}
     \end{pmatrix}
        }
        {\eede'+\frac{\hhdeov}{\eedeov}\hhdet'+\frac{\phhdeov}{\eedeov}\phhdet'} &
\\
 + \frac{\kappa\tov}{\eedeov}\frac{\vmrov}{\pvmrov} 
   \frac{\frac{1}{\eede'+\frac{\hhdeov}{\eedeov}\hhdet'+
                            \frac{\phhdeov}{\eedeov}\phhdet'}}
     {1+\frac{\kappa\eedeov}{c}(\eede'+\frac{\hhdeov}{\eedeov}\hhdet'+
                                 \frac{\phhdeov}{\eedeov}\phhdet')}  \vmrt' 
+\frac{\kappa\tov}{\eedeov}\frac{\frac{1}{\eede'+\frac{\hhdeov}{\eedeov}\hhdet'+
                            \frac{\phhdeov}{\eedeov}\phhdet'}}
     {1+\frac{\kappa\eedeov}{c}(\eede'+\frac{\hhdeov}{\eedeov}\hhdet'+
                                 \frac{\phhdeov}{\eedeov}\phhdet')}  \pvmrt'& =
\\
   \frac{\mu\tov}{\eedeov}
   \frac{\frac{1}{\eede'+\frac{\hhdeov}{\eedeov}\hhdet'+\frac{\phhdeov}{\eedeov}\phhdet'}}
        {1+\frac{\mu\eedeov}{c}
        (\eede'+\frac{\hhdeov}{\eedeov}\hhdet'+\frac{\phhdeov}{\eedeov}\phhdet')} 
        (\frac{\vwnov}{\pvmrov}  \vwnt'- \frac{\vmrov}{\pvmrov}\vmrt')- 
   \frac{\mu\tov}{\eedeov}
   \frac{\frac{1}{\eede'+\frac{\hhdeov}{\eedeov}\hhdet'+\frac{\phhdeov}{\eedeov}\phhdet'}}
        {1+\frac{\mu\eedeov}{c}
        (\eede'+\frac{\hhdeov}{\eedeov}\hhdet'+\frac{\phhdeov}{\eedeov}\phhdet')} \pvmrt'
    + \frac{\fxtov\tov}{\pvmrov}\fxt' .&
\end{aligned}
\label{ppmtl4}
\end{align}

\subsection{Parameter size and rescaled equations}
In this subsection, we fix the characteristic values. As set out in 
the previous subsection, we shall choose mean values or maximal values
of the concerned physical quantities.
The parameter  $\tov$ is the observation time scale. We consider
that it is about several months. Then we set 
\begin{equation}
\tov \sim 100 \text{ days } \sim 2400 h,
\end{equation}
beside this $\omov$ is the tide frequency, meaning $1/\omov$ 
is the tide duration, i.e.:
\begin{equation}\label{PS1}
\frac{1}{\omov} \sim 13 h.
\end{equation}
Hence we exhibit a small parameter:
\begin{equation}\label{PS2}
\eps = \frac{1}{\tov \omov} \sim  \frac{1}{200}.
\end{equation}

Then, we make a strong assumption, which is that 
$\pvmrt$ and  $\phhdet$ are really perturbations. In other words, we consider
that
\begin{equation}
\frac{\pvmrov}{\vmrov} \sim \frac{\phhdeov}{\hhdeov} \sim  \eps.
\end{equation}

~

Concerning the Coriolis parameter, for latitudes about $45^\circ$
$f\sim \pi/day \sim 4\; 10^{-5} /s$,
then $f\tov\sim \pi/2\eps$.
Concerning the other parameters of physical meaning, several choices are
possible, according to the turbulence action, the nature of the ocean bottom, 
the shape of the ocean free surface and so on. We focus on one of those choices, 
being aware that others, that would lead to other models, are also reasonable.
For the viscosity, we chose the value of the water viscosity at $20^\circ C$, i.e.
$c\sim10^{-2} cm^2/s \sim 10^{-7} km^2/day$, then
$c\tov\sim 10^{-5} km^2$. 
Concerning  the friction coefficients, the bottom friction coefficient is
$\kappa\sim10^{-4}m/s \sim 10^{-2} km/day$ and the air-water friction coefficient
is $\mu\sim10^{-6}m/s \sim10^{-4} km/day$.
Those values are consistent with the ones used, for instance, in
Dawson and Proft \cite{DawPro} or
Gerbeau and Perthame \cite{GerPer}.
Then $\kappa\tov\sim 1 km$, $\mu\tov\sim 10^{-2} km$,
$\kappa/c\sim10^{5}/km$ and $\mu/c\sim10^{3}/km$.
We also have $g\sim 10 m/s^2\sim 10^{6} km/day^2$ and $g\tov\sim 10^{8} km/day$.

~

We now turn to the ratios determining the asymptotic analysis we
have to realize.
Having in mind our final target, i.e. the drift of things in the ocean over
long time periods, we notice that such a drift may take place relatively 
far from the coast, above the  \textbf{continental shelf}. It may also
take place in a large and relatively closed bay, with a long residence time
of the ocean water. Such a domain will be called \textbf{coastal zone}.
As was the case in 1999 / 2000 for the Erika oil slick along
the French Atlantic coast, the drift may occur for weeks along a thin
\textbf{layer} following the coast. Those remarks guide the choices concerning
the geometrical assumptions we consider further.  

$\llov$ and  $\lov$ represent the characteristic lengths of the domain
where the drift takes place and $\vmrov/\omov$ the characteristic
distance the water covers in the tide duration. 
Following 
Salomon and Breton \cite{SaBre},
\u{C}etina, Rajar and Povinec \cite{CeRaPo},
Bao, Gao and Yan \cite{BaoGaoYan} or
Cai, Huang and Long \cite{CaiHuLon},
we can state that
this distance is about a few kilometers (from $5$ to $20$) in the cases we are interested in.
If the domain under consideration is a \textbf{continental shelf}, 
the characteristic sea water velocity $\vmrov$ is about $0.5 km/h$ and 
then, we have
\begin{equation}
\frac{\vmrov}{\omov} \sim 5 km,
\end{equation}
and if we set 
$\llov\sim$ $\lov\sim$ $500 km$, $\eedeov\sim 300 m$ and $\hhdeov\sim 3m$, then
\begin{equation}
\frac{\frac{\vmrov}{\omov}}{\llov} \sim 2 \eps ,
~~~~~~ \frac{\hhdeov}{\eedeov} \sim 2 \eps.
\end{equation}
We also have  
$\pvmrov\sim \eps \vmrov\sim 2.5\;10^{-3} km/h \sim 6\;10^{-2} km/day$,
then we get $g\tov/\pvmrov \sim 1.7\; 10^{9}$.
Since $\phhdeov\sim \eps\hhdeov\sim 1.5\; 10^{-2} m$,
we obtain $\phhdeov/\llov\sim 3\;10^{-8}$.
Hence  
\begin{equation}\label{REr2}
\frac{g\tov}{\pvmrov}\frac{\phhdeov}{\llov}\sim 50 \sim \frac{1}{4\eps}.
\end{equation}
Moreover 
\begin{gather}
\frac{c\tov}{\llov^2}\sim \frac{10^{-5}}{25\;10^{4} }\sim 13\eps^5,\label{REr3}
\\
\frac{\kappa\tov}{\eedeov} \sim \frac{1}{0.3}\sim 3.3\sim 3,
~~~~~
\frac{\kappa\eedeov}{c} \sim 3\; 10^{4} \sim  \frac{0.8}{\eps^2}
\\
\frac{\mu\tov}{\eedeov} \sim \frac{10^{-2}}{0.3}\sim 3.3\;10^{-2}\sim 6\eps,
~~~~~
\frac{\mu\eedeov}{c} \sim 3\; 10^{2}  \sim \frac{1.5}{\eps}.\label{REr5}
\end{gather}
Concerning the wind velocity, when the weather is calm, $10 km/h$ is a relevant
characteristic value, while $100 km/h$ may be a good choice in 
stormy conditions. Hence, we shall consider 
\begin{equation}
\frac{\vmrov}{\vwnov}\sim \frac{0.5}{10}\sim \frac{1}{20},
\end{equation}
in calm weather regime, and, 
\begin{equation}
\frac{\vmrov}{\vwnov}\sim \frac{0.5}{100}\sim \eps,
\end{equation}
in storm regime.
Expressing now the following ratios
\begin{equation}\label{REr1}
\frac{\vmrov\tov}{\llov}\sim
\tov \omov\frac{\frac{\vmrov}{\omov}}{\llov}, ~~~~~~~
\frac{\pvmrov\tov}{\llov}\sim \frac{\vmrov\tov}{\llov}
\frac{\pvmrov}{\vmrov},
\end{equation}
and moreover setting $\fxtov\tov\sim\pvmrov$ and removing the ' for clarity,
we can write the rescaled equation (\ref{ppmtl3})-(\ref{ppmtl4}) in
the case of a continental shelf:
\begin{align} \ds
&
\begin{aligned}
  \fracp{\phhdet}{\tc} + 
  \nabla (\frac{1}{\eps}\eede+2\hhdet))\cdot\pvmrt+
  (\frac{1}{\eps}\eede+2\hhdet)\nabla\cdot\pvmrt+
  2(\nabla \phhdet)\cdot\vmrt +
  2\phhdet (\nabla\cdot\vmrt) &
\\ \ds + 
  2\eps\big( (\nabla \phhdet)\cdot\pvmrt+
  \phhdet(\nabla\cdot\pvmrt)\big)&=0, 
\end{aligned}
\label{REcs3}
\\ \ds
&
\begin{aligned}
\fracp{\pvmrt}{\tc} +
  2(\nabla\pvmrt)\vmrt+
  2(\nabla\vmrt)\pvmrt+
  2\eps (\nabla\pvmrt)\pvmrt+
  \frac{\pi}{2\eps}\pvmrt^\perp+
  \frac{1}{4\eps}\nabla\phhdet-
  13\eps^4\Delta  \vmrt-
  13\eps^5 \Delta  \pvmrt ~~~~&
\\
  -13\eps^4
  \frac{(\nabla \vmrt)\nabla(\eede+2\eps\hhdet)}
       {\eede+2\eps\hhdet +2\eps^2\phhdet}
 -26\eps^6
     \frac{(\nabla\vmrt)\nabla\phhdet }
          {\eede+2\eps\hhdet+2\eps^2\phhdet}
 -13\eps^5
   \frac{(\nabla\pvmrt ) \nabla(\eede+2\eps\hhdet)}
        {\eede+2\eps\hhdet+2\eps^2\phhdet}&
\\
 -26\eps^7
   \frac{(\nabla\pvmrt ) \nabla\phhdet}
        {\eede+2\eps\hhdet+2\eps^2\phhdet}&
\\
+\frac3\eps \frac{\frac{1}{\eede+2\eps\hhdet+ 2\eps^2\phhdet}}
     {1+\frac{0.8}{\eps^2}(\eede+2\eps\hhdet+2\eps^2\phhdet)} \vmrt
+3\frac{\frac{1}{\eede+2\eps\hhdet+ 2\eps^2\phhdet}}
     {1+\frac{0.8}{\eps^2}(\eede+2\eps\hhdet+2\eps^2\phhdet)}  \pvmrt
 &= 
\\ 
   6\frac{\frac{1}{\eede+2\eps\hhdet+2\eps^2\phhdet}}
        {1+\frac{1.5}{\eps}
        (\eede+2\eps\hhdet+2\eps^2\phhdet)} 
        (\gamma\vwnt- \vmrt)- 
   6\eps\frac{\frac{1}{\eede+2\eps\hhdet+2\eps^2\phhdet}}
        {1+\frac{1.5}{\eps}
        (\eede+2\eps\hhdet+2\eps^2\phhdet)} \pvmrt
    &+\fxt .
\end{aligned}
\label{REcs4}
\end{align}
where $\gamma=20=1/(10\eps)$ in calm weather regime.
In storm regime, which is what it is supposed for model 
(\ref{REcs3intr})-(\ref{REcs4intr}) presented 
in the introduction, $\gamma=1/\eps$.

~

If the domain is a \textbf{coastal zone},  $\vmrov \sim 1 km/h$
\begin{equation}
\frac{\vmrov}{\omov} \sim 10 km,
\end{equation} 
and we set 
$\llov\sim$ $\lov\sim$ $5 km$, $\eedeov\sim 50 m$ et $\hhdeov\sim 10m$. In this case
\begin{equation}
\frac{\frac{\vmrov}{\omov}}{\llov} \sim 2 ,~~~~~~ \frac{\hhdeov}{\eedeov} \sim \frac15.
\end{equation}

We also have  
$\pvmrov\sim \eps \vmrov\sim 1.2\;10^{-1} km/day$,
then we get $g\tov/\pvmrov \sim 8\; 10^{8}$.
Since $\phhdeov\sim \eps\hhdeov\sim 5\; 10^{-2} m$
we obtain $\phhdeov/\llov\sim 10^{-5}$.
Hence 
\begin{equation}
\frac{g\tov}{\pvmrov}\frac{\phhdeov}{\llov}\sim 8\;10^{3}
\sim  \frac{0.2}{\eps^2}.
\end{equation}
Moreover 
\begin{gather}
\frac{c\tov}{\llov^2}\sim \frac{10^{-5}}{25}\sim 0.6\eps^3,
\\
\frac{\kappa\tov}{\eedeov} \sim \frac{1}{0.05}\sim 20\sim \frac{1}{10\eps},
~~~~~
\frac{\kappa\eedeov}{c} \sim 5\; 10^{3} \sim  \frac{1}{10\eps^2},
\\
\frac{\mu\tov}{\eedeov} \sim \frac{10^{-2}}{0.05}\sim 0.2,
~~~~~
\frac{\mu\eedeov}{c} \sim 50  \sim \frac{1}{4\eps}.
\end{gather}
Concerning the wind velocity, we have 
\begin{equation}
\frac{\vmrov}{\vwnov}\sim \frac{1}{10},
\end{equation}
in calm weather regime, and, in storm regime
\begin{equation}
\frac{\vmrov}{\vwnov}\sim \frac{1}{100}\sim 2\eps.
\end{equation}
Hence, the rescaled equation reads in this case:
\begin{align} \ds
&
\begin{aligned}
  \fracp{\phhdet}{\tc} + 
  \frac{2}{\eps}(\nabla (5\eede+\hhdet))\cdot\pvmrt+
  \frac{2}{\eps}(5\eede+\hhdet)\nabla\cdot\pvmrt+
  \frac{2}{\eps}(\nabla \phhdet)\cdot\vmrt +
  \frac{2}{\eps}\phhdet (\nabla\cdot\vmrt) &
\\ \ds + 
   2(\nabla \phhdet)\cdot\pvmrt+
  2\phhdet(\nabla\cdot\pvmrt)&=0, 
\end{aligned}
\label{REcz3}
\\ \ds
&
\begin{aligned}
\fracp{\pvmrt}{\tc} +
  \frac{2}{\eps}(\nabla\pvmrt)\vmrt+
  \frac{2}{\eps}(\nabla\vmrt)\pvmrt+
  2(\nabla\pvmrt)\pvmrt+
  \frac{\pi}{2\eps}\pvmrt^\perp+
   \frac{0.2}{\eps^2}\nabla\phhdet-
  0.6\eps^2 \Delta  \vmrt -
   0.6\eps^3 \Delta  \pvmrt&
\\
  -0.6\eps^2
  \frac{(\nabla \vmrt)\nabla(\eede+\hhdet) }
       {\eede+\frac{1}{5}\hhdet +\frac{\eps}{5}\phhdet}
 -0.1\eps^2
     \frac{(\nabla\vmrt)\nabla\phhdet }
          {\eede+\frac{1}{5}\hhdet+\frac{\eps}{5}\phhdet'}
 -0.6\eps^2
   \frac{(\nabla\pvmrt ) \nabla(\eede+\hhdet)}
        {\eede+\frac{1}{5}\hhdet+\frac{\eps}{5}\phhdet}&
\\
 -0.1\eps^3
   \frac{(\nabla\pvmrt ) \nabla\phhdet}
        {\eede+\frac{1}{5}\hhdet+\frac{\eps}{5}\phhdet}&
\\+\frac{1}{10}\frac{\frac{1}{\eede+\frac{1}{5}\hhdet+ \frac{\eps}{5}\phhdet}}
     {1+\frac{1}{10\eps^2}(\eede+\frac{1}{5}\hhdet+\frac{\eps}{5}\phhdet)}  \vmrt
+\frac{1}{10\eps}\frac{\frac{1}{\eede+\frac{1}{5}\hhdet+ \frac{\eps}{5}\phhdet}}
     {1+\frac{1}{10\eps^2}(\eede+\frac{1}{5}\hhdet+\frac{\eps}{5}\phhdet)}  \pvmrt
 &= 
\\ 
   0.2\frac{\frac{1}{\eede+\frac{1}{5}\hhdet+\frac{\eps}{5}\phhdet}}
        {1+\frac{1}{4\eps}
        (\eede+\frac{1}{5}\hhdet+\frac{\eps}{5}\phhdet)} 
        (\frac{\gamma}{2\eps} \vwnt - \unsep\vmrt)- 
   0.2\frac{\frac{1}{\eede+\frac{1}{5}\hhdet+\frac{\eps}{5}\phhdet}}
        {1+\unsep
        (\eede+\frac{1}{5}\hhdet+\frac{\eps}{5}\eps\phhdet)} \pvmrt
    &+\fxt .
\end{aligned}
\label{REcz4}
\end{align}
where $\gamma/2=10=1/20\eps$ in still weather and $\gamma/2=1/2\eps$
in stormy weather.

~

We will give the name of \textbf{coastal layer} to a domain
having the following characteristics
\begin{equation}
\frac{\vmrov}{\omov} \sim 10 km,
\end{equation} 
and
$\llov\sim$ $500 km$, $\lov\sim$ $5 km$, $\eedeov\sim 50 m$ et $\hhdeov\sim 10m$ and then  
\begin{equation}
\frac{\frac{\vmrov}{\omov}}{\llov} \sim 4 \eps ,~~~~~~ 
\frac{\lov}{\llov} \sim 2 \eps,~~~~~~ 
\frac{\hhdeov}{\eedeov} \sim \frac15 .
\end{equation} 
In this case, we have 
$\pvmrov\sim \eps \vmrov\sim 1.2\;10^{-1} km/day$,
then we get $g\tov/\pvmrov \sim 8\; 10^{8}$.
Since $\phhdeov\sim \eps\hhdeov\sim 5\; 10^{-2} m$,
we obtain $\phhdeov/\llov\sim 10^{-7}$.
Hence 
\begin{equation}
\frac{g\tov}{\pvmrov}\frac{\phhdeov}{\llov}\sim 80
\sim  \frac{0.4}{\eps}.
\end{equation}
Moreover 
\begin{gather}
\frac{c\tov}{\llov^2}\sim \frac{10^{-5}}{25\;10^{4} }\sim 13\eps^5,
\\
\frac{\kappa\tov}{\eedeov} \sim \frac{1}{0.05}\sim 20\sim \frac{1}{10\eps},
~~~~~
\frac{\kappa\eedeov}{c} \sim 5\; 10^{3} \sim  \frac{1}{10\eps^2},
\\
\frac{\mu\tov}{\eedeov} \sim \frac{10^{-2}}{0.05}\sim 0.2,
~~~~~
\frac{\mu\eedeov}{c} \sim 50  \sim \frac{1}{4\eps}.
\end{gather}
The considerations concerning the wind velocity are the same as in the case
of a coastal zone. The the rescaled equation for a coastal layer writes: 
\begin{align} \ds
&
\begin{aligned}
 \fracp{\phhdet}{\tc} + 
  \begin{pmatrix}\ds 20\fracp{\eede}{\xc_1}+
    \fracp{\hhdet}{\xc_1}\\ \ds
    \frac{2}{\eps}(5\fracp{\eede}{\xc_2}+\fracp{\hhdet}{\xc_2})
  \end{pmatrix}
  \cdot\pvmrt+
  (20\eede+\hhdet)
  \big(\fracp{\pvmrt_1}{\xc_1}+\frac{1}{2\eps}\fracp{\pvmrt_2}{\xc_2}\big) + 
  4\begin{pmatrix}\ds\fracp{\phhdet}{\xc_1}\\ 
         \ds\frac{1}{2\eps}\fracp{\phhdet}{\xc_2}
  \end{pmatrix}
  \cdot\vmrt~~~~~~& \\
  +4\phhdet\big(\fracp{\vmrt_1}{\xc_1}+\frac{1}{2\eps}\fracp{\vmrt_2}{\xc_2}\big)
  + \begin{pmatrix}\ds4\eps\fracp{\phhdet}{\xc_1}\\
                             \ds 2\fracp{\phhdet}{\xc_2}\end{pmatrix}
  \cdot\pvmrt+  
  \phhdet\big(4\eps\fracp{\pvmrt_1}{\xc_1}+2\fracp{\pvmrt_2}{\xc_2}\big) &=0, 
\end{aligned}
\label{REcl3}
\end{align}
%\displaybreak[5]
%\\
\begin{align}
&\begin{aligned}
  \fracp{\pvmrt}{\tc} +
  \Big( 4\fracp{\pvmrt}{\xc_1},\frac{2}{\eps}\fracp{\pvmrt}{\xc_2}\Big)\vmrt+
  \Big( 4\fracp{\vmrt}{\xc_1},\frac{2}{\eps}\fracp{\vmrt}{\xc_2}\Big)\pvmrt+
     \Big( 4\eps\fracp{\pvmrt}{\xc_1},2\fracp{\pvmrt}{\xc_2}\Big)\pvmrt
  +\frac{\pi}{2\eps}\pvmrt^\perp~~~~~~~~~~~~~~~~~~~~~&
\\
  + \frac{0.4}{\eps}\begin{pmatrix}
     \ds\fracp{\phhdet}{\xc_1}\\\ds\frac{1}{2\eps}\fracp{\phhdet}{\xc_2}
   \end{pmatrix}
  - (13\eps^4\fracpse{\vmrt}{{\xc_1}}+\frac{13\eps^2}{4}\fracpse{\vmrt}{{\xc_2}})
  - (13\eps^5\fracpse{\pvmrt}{{\xc_1}}+\frac{13\eps^3}{4}\fracpse{\pvmrt}{{\xc_2}})&
\\
   - 13\eps^4
  \frac{ \ds
  \Big( \fracp{\vmrt}{\xc_1},\frac{1}{2\eps}\fracp{\vmrt}{\xc_2}\Big)
  \begin{pmatrix}\ds\fracp{\eede}{\xc_1}+
    \frac 15 \fracp{\hhdet}{\xc_1}\\ \ds
    \frac{1}{2\eps}(\fracp{\eede}{\xc_2}+
    \frac 15 \fracp{\hhdet}{\xc_2})
  \end{pmatrix}
       }
       {\eede+\frac 15 \hhdet +\frac \eps 5 \phhdet}&
\\ 
   -\frac{13\eps^5}{5}
     \frac{\ds 
    \Big( \fracp{\vmrt}{\xc_1},\frac{1}{2\eps}\fracp{\vmrt}{\xc_2}\Big)
    \begin{pmatrix}
     \ds\fracp{\phhdet}{\xc_1}\\\ds\frac{1}{2\eps}\fracp{\phhdet}{\xc_2}
   \end{pmatrix}
          }
          {\eede+\frac 15 \hhdet+\frac{\eps}{5}\phhdet}
   -13 \eps^5
   \frac{ \ds
     \Big( \fracp{\pvmrt}{\xc_1},\frac{1}{2\eps}\fracp{\pvmrt}{\xc_2}\Big)
     \begin{pmatrix}\ds\fracp{\eede}{\xc_1}+
      \frac 15 \fracp{\hhdet}{\xc_1}\\ \ds
      \frac{1}{2\eps}(\fracp{\eede}{\xc_2}+
      \frac 15\fracp{\hhdet}{\xc_2})
     \end{pmatrix}
        }
        {\eede+\frac 15\hhdet+\frac{\eps}{5}\phhdet}&
\\
   -\frac{13\eps^6}{5}
   \frac{ \ds
     \Big( \fracp{\pvmrt}{\xc_1},\frac{1}{2\eps}\fracp{\pvmrt}{\xc_2}\Big)
     \begin{pmatrix}
       \ds\fracp{\phhdet}{\xc_1}\\\ds\frac{1}{2\eps}\fracp{\phhdet}{\xc'_2}
     \end{pmatrix}
        }
        {\eede+\frac 15\hhdet+\frac{\eps}{5}\phhdet} &
\\
 +\frac{1}{10\eps^2}\frac{\frac{1}{\eede+\frac 15\hhdet+ \frac{\eps}{5}\phhdet}}
     {1+\frac{1}{10\eps^2}(\eede+\frac 15\hhdet+\frac{\eps}{5}\phhdet)} \vmrt
+\frac{1}{10\eps}\frac{\frac{1}{\eede+\frac 15\hhdet+ \frac{\eps}{5}\phhdet}}
     {1+\frac{1}{10\eps^2}(\eede+\frac 15\hhdet+\frac{\eps}{5}\phhdet)}  \pvmrt& =
\\
   0.2\frac{\frac{1}{\eede+\frac 15\hhdet+\frac{\eps}{5}\phhdet}}
        {1+\frac{1}{4\eps}
        (\eede+\frac 15\hhdet+\frac{\eps}{5}\phhdet)} 
        (\frac{\gamma}{2\eps}\vwnt- \unsep\vmrt)- 
   0.2\frac{\frac{1}{\eede+\frac 15\hhdet+\frac{\eps}{5}\phhdet}}
        {1+\frac{1}{4\eps}
        (\eede+\frac 15\hhdet+\frac{\eps}{5}\phhdet)} \pvmrt
    +&\fxt .
\end{aligned}
\label{SM1cl4}
\end{align}

\section{Existence}
\label{ExisSec}
\subsection{Simplified system for continental shelf}

In this section, we focus on one of the models introduced in the previous section, 
and we explore some of its mathematical properties. 

More precisely, we consider a simplified version of 
system (\ref{REcs3})-(\ref{REcs4})
which consists in considering that the ocean bottom is flat, i.e. 
$\eede\equiv 1$, in forgetting all the power of $\eps$ greater than $1$
and in setting all constants to 1. Then we obtain system 
(\ref{SM1cs3intr})-(\ref{SM1cs2intr}) and we prove an existence result
for the solution of this system.

Although the results given in sections \ref{ExisSec} and \ref{AsBe} 
are specific to the model (\ref{SM1cs3intr})-(\ref{SM1cs2intr}), 
we expect that similar methods could be used to prove similar results
for the other models introduced in section \ref{ModSec}.

%\begin{align} \ds
%&
%\begin{aligned}
%  \fracp{\phhdet}{\tc} + 
%  (\nabla \hhde)\cdot\pvmrt+
%  (\frac{1}{\eps}+\hhde)(\nabla\cdot\pvmrt)+
%  (\nabla \phhdet)\cdot\vmr +
%  \phhdet (\nabla\cdot\vmr) + 
%  \eps\big( (\nabla \phhdet)\cdot\pvmrt+
%  \phhdet(\nabla\cdot\pvmrt)\big)&=0, 
%\end{aligned}
%\label{SM1cs3}
%\\ \ds
%&
%\fracp{\pvmrt}{\tc} +
%  (\nabla\pvmrt)\vmr+
%  (\nabla\vmr)\pvmrt+
%  \eps (\nabla\pvmrt)\pvmrt+
%  \frac{1}{\eps}\pvmrt^\perp+
%  \frac{1}{\eps}\nabla\phhdet
% =\vwn.
%\label{SM1cs2}
%\end{align}
% 
\subsection{Proof of Theorem \ref{THM1}}
Setting $\uvec=(\phhdet,\pvmrt)=(\phhdet,\pvmrt_1,\pvmrt_2)$, $\uvec^\perp=(0,\pvmrt^\perp)$
and introducing
\begin{gather}
B^1(t,\tsep,\xvec,\eps \uvec)
=\begin{pmatrix}
\vmrt_1+\eps\pvmrt_1 &\ds \frac{1}{\eps}+\hhdet+\eps\phhdet & 0 \\ \ds
\frac{1}{\eps} & \vmrt_1+\eps\pvmrt_1 & 0 \\
0 & 0 &\ds \vmrt_1+\eps\pvmrt_1
 \end{pmatrix},
\\
B^2(t,\tsep,\xvec,\eps \uvec)
=\begin{pmatrix}
\vmrt_2+\eps\pvmrt_2 & 0 &\ds \frac{1}{\eps}+\hhdet+\eps\phhdet \\
0 & \vmrt_2+\eps\pvmrt_2 & 0 \\ \ds
\frac{1}{\eps}& 0 & \vmrt_2+\eps\pvmrt_2
 \end{pmatrix},
\end{gather}
and
\begin{gather}
F(t,\tsep,\xvec, \uvec)
=\begin{pmatrix}\ds
-(\fracp{\hhdet}{x_1} \pvmrt_1 + \fracp{\hhdet}{x_2} \pvmrt_2 )
-(\fracp{\vmrt_1}{x_1}+\fracp{\vmrt_2}{x_2})\phhdet
\\ \ds
\vwnt_1 - (\fracp{\vmrt_1}{x_1}\pvmrt_1+\fracp{\vmrt_1}{x_2}\pvmrt_2)
\\ \ds
\vwnt_2 - (\fracp{\vmrt_2}{x_1}\pvmrt_1+\fracp{\vmrt_2}{x_2}\pvmrt_2)
 \end{pmatrix},
\end{gather}
equations (\ref{SM1cs3intr})-(\ref{SM1cs2intr}) read
\begin{gather}
\label{PresqueSystHyp}
\fracp{\uvec}{t} + B^1 \fracp{\uvec}{x_1} + B^2 \fracp{\uvec}{x_2} + \frac{1}{\eps}\uvec^\perp= F.
\end{gather}
And, introducing 
\begin{gather}
\label{A0def}
A^0(t,\tsep,\xvec,\eps^2\phhdet) = 
\begin{pmatrix}
\ds \frac{1}{1+\eps\hhdet+\eps^2\phhdet} & 0 & 0 \\
0 & 1 & 0 \\
0 & 0 & 1
\end{pmatrix},
\end{gather}
\begin{gather}
A^1(t,\tsep,\xvec,\eps \uvec)
=\begin{pmatrix}\ds
\frac{\vmrt_1+\eps\pvmrt_1}{1+\eps\hhdet+\eps^2\phhdet} & 0 & 0 \\ \ds
0 & \vmrt_1+\eps\pvmrt_1 & 0 \\
0 & 0 &\ds \vmrt_1+\eps\pvmrt_1
 \end{pmatrix},
\end{gather}
\begin{gather}
A^2(t,\tsep,\xvec,\eps \uvec)
=\begin{pmatrix}\ds
\frac{\vmrt_2+\eps\pvmrt_2}{1+\eps\hhdet+\eps^2\phhdet} & 0 & 0\\
0 & \vmrt_2+\eps\pvmrt_2 & 0 \\ \ds
0 & 0 & \vmrt_2+\eps\pvmrt_2
 \end{pmatrix},
\end{gather}
\begin{gather}
S^1=\begin{pmatrix}
0 & 1 & 0 \\
1 & 0 & 0 \\
0 & 0 & 0
    \end{pmatrix} \text{ ~ and ~ }
S^2=\begin{pmatrix}
0 & 0 & 1 \\
0 & 0 & 0 \\
1 & 0 & 0
    \end{pmatrix},
\end{gather}
equation (\ref{PresqueSystHyp}) yields the following symmetric hyperbolic system:
\begin{gather}
\label{SystHyp}
A^0\fracp{\uvec}{t} 
+ A^1 \fracp{\uvec}{x_1} 
+ A^2 \fracp{\uvec}{x_2}
+ \frac{1}{\eps} S^1 \fracp{\uvec}{x_1} 
+ \frac{1}{\eps} S^2 \fracp{\uvec}{x_2}
+ \frac{1}{\eps}\uvec^\perp
=F_0 = A^0 F.
\end{gather}
Hence applying Kato \cite{Kato:1975} or Majda \cite{Majda:1984},
we deduce that for any $\eps$, under the assumptions of Theorem \ref{THM1}, 
the classical solution of (\ref{SM1cs3intr}), (\ref{SM1cs2intr}) and (\ref{Cini}) 
exists and is unique on a time interval. what remains to prove is that this time
interval does not depend on $\eps$.

~

Before going further in the proof, we may observe the following differences 
between equation (\ref{SystHyp}) and 
the type of problems, also depending on a small parameter, studied in 
Klainerman and Majda \cite{KlaiMaj:1981,KlaiMaj:1982},
Schochet \cite{schochet:1986, schochet11:1986, schochet:1988},
and Metivier and Schochet \cite{MetScho}.
Some of those differences simplify the problem: the nonlinearity in $A^0$ are 
functions of only $\eps^2 \phhdet$ and in $A^1$ and $A^2$ of only $\eps \uvec$.
Some others make the results proved by those authors unable to be applied directly:  
$A^0$, $A^1$ and $A^2$ depend on
$t/\eps$ and the singular term $\uvec^\perp/\eps$  involves the function $\uvec$ itself
and not order 1 derivatives of it. 
Nonetheless, the now classical calculus procedures carried out 
in the concerned papers may be followed
in order to obtain the right estimates allowing for the conclusion. We sketch the concerned
computations hereafter.

We set $\alpha = (\alpha_1,\alpha_2) \in \nit^2$ with 
$|\alpha|=\alpha_1+\alpha_2 \leq s$
and $\ds D^\alpha \uvec = \frac{\partial^{|\alpha|} \uvec}
{{\partial x_1}^{\alpha_1}{\partial x_2 }^{\alpha_2}}$.
Applying $D^\alpha$ to equation (\ref{SystHyp}) yields
\begin{gather}
\label{DalSystHyp}
A^0\fracp{D^\alpha\uvec}{t} 
+ A^1 \fracp{D^\alpha\uvec}{x_1} 
+ A^2 \fracp{D^\alpha\uvec}{x_2}
+ \frac{1}{\eps} S^1 \fracp{D^\alpha\uvec}{x_1} 
+ \frac{1}{\eps} S^2 \fracp{D^\alpha\uvec}{x_2}
+ \frac{1}{\eps}(D^\alpha\uvec)^\perp
= F^\alpha,
\end{gather}
with
\begin{gather}
\label{Falpha}
F^\alpha = D^\alpha  F_0
- [D^\alpha , A^0\fracp{}{t} ]\uvec 
- [D^\alpha , A^1\fracp{}{x_1} ]\uvec 
- [D^\alpha , A^2\fracp{}{x_2}]\uvec,
\end{gather}
$[~, ~]$ standing for the classical commutator.

Multiplying equation (\ref{DalSystHyp}) by $2 D^\alpha\uvec$, 
, integrating on $\rit^2$ and noticing that
\begin{gather}
2 \int A^0\fracp{D^\alpha\uvec}{t}\cdot D^\alpha\uvec \,d\xvec 
= \frac{\ds d \Big(\int  A^0 D^\alpha\uvec \cdot D^\alpha\uvec \,d\xvec \Big)}{dt}
- \int  \frac{d(A^0)}{dt} D^\alpha\uvec \cdot D^\alpha \uvec \,d\xvec ,
\\
2 \int A^j\fracp{D^\alpha\uvec}{x_j}\cdot D^\alpha\uvec \,d\xvec
= \ds \Big(\int  \frac{\ds d \big( A^j D^\alpha\uvec \cdot D^\alpha\uvec \big) }{dx_j} \,d\xvec \Big)
- \int  \frac{d(A^j)}{dx_j} D^\alpha\uvec \cdot D^\alpha\uvec \,d\xvec 
\nonumber
\\
~~~~~~~~~~~~~~~~~~~~~~~~~~~~~~~~~~~~~~~~~~~~~~~~~~~~~~~~~
= - \int  \frac{d(A^j)}{dx_j} D^\alpha\uvec \cdot D^\alpha\uvec \,d\xvec 
\\
2 \int S^j\fracp{D^\alpha\uvec}{x_j}\cdot D^\alpha\uvec \,d\xvec
= -2 \int S^j\fracp{D^\alpha\uvec}{x_j}\cdot D^\alpha \,d\xvec = 0,
\end{gather}
for $j=1,2$, and
\begin{gather}
2 \int (D^\alpha\uvec)^\perp \cdot D^\alpha\uvec \,d\xvec = 0,
\end{gather}
we obtain
\begin{multline}
\label{est1}
\frac{\ds d \Big(\int  A^0 D^\alpha\uvec \cdot D^\alpha\uvec \,d\xvec \Big)}{dt}
= \int  \frac{d(A^0)}{dt} D^\alpha\uvec \cdot D^\alpha\uvec \,d\xvec 
\\
+ \int  \frac{d(A^1)}{dx_1} D^\alpha\uvec \cdot D^\alpha\uvec \,d\xvec
+ \int  \frac{d(A^2)}{dx_2} D^\alpha\uvec \cdot D^\alpha\uvec \,d\xvec
+ \int F^\alpha \cdot D^\alpha\uvec \,d\xvec.
\end{multline}
For all the estimates to come, all the constants which are needed are called $c$. 
Since the dependency of $A^0$ with respect to $t/\eps$ is done through $\eps\hhdet$,
and since $s>3$,
we can deduce that for any $t$ and $\xvec$,
\begin{gather}
\label{EstIneg1}
\Big| \frac{d(A^0_{11}(t,\tsep,\xvec,\eps^2\phhdet))}{dt} \Big| 
\leq c (1+ \eps^2 \big | \fracp{\phhdet}{t} \big | )
\leq c (1+ \eps^2 \sup_{\xvec\in\rit^2}\big | \fracp{\phhdet}{t} \big | )
\leq c (1+ \eps^2 \big \| \fracp{\phhdet}{t} \big \|_{s-1} ),
\end{gather}
where $\| ~\|_{s-1}$ stands for the norm in $H^{s-1}(\rit^2)$.
The time derivatives of the other entries of $A^0$ are zero.
Hence the first term of the right hand side of (\ref{est1}) may be estimated
\begin{gather}
\label{EstIneg2}
\bigg|\int  \frac{d(A^0)}{dt} D^\alpha\uvec \cdot D^\alpha\uvec \,d\xvec \bigg|
\leq c (1+ \eps^2 \big \| \fracp{\phhdet}{t} \big \|_{s-1} ) \,
\| D^\alpha\uvec \|_{0}^2 
\leq c (1+ \eps^2 \big \| \fracp{\phhdet}{t} \big \|_{s-1} ) \,
\|  \uvec \|_{s}^2,
\end{gather}
where $\| ~\|_{s}$ stands for the norm in $(H^{s}(\rit^2))^2$
and $\| ~\|_{0}$ for the norm in $(L^2(\rit^2))^2$.
Concerning the entries $A^i_{kl}$ of $A^i$ for $i=1,2$, 
\begin{multline}
\label{EstIneg3}
\Big| \frac{d(A^i_{kl}(t,\tsep,\xvec,\eps \uvec))}{dx_i} \Big| 
\leq c (1+ \eps \big | \fracp{\pvmrt_i}{x_i}\big | +\eps^2 \big | \fracp{\phhdet}{x_i} \big | )
\leq c (1+ \eps \sup_{\xvec\in\rit^2}\big | \fracp{\pvmrt_i}{x_i}\big | 
+ \eps^2 \sup_{\xvec\in\rit^2}\big | \fracp{\phhdet}{x_i} \big | )
\\
\leq c (1+ \eps \big \| \fracp{\uvec }{x_i} \big \|_{s-1} )\,
\leq c (1+ \eps \big \| \uvec  \big \|_{s} ).
\end{multline}
Hence 
\begin{gather}
\label{EstIneg4}
\bigg|\int  \frac{d(A^i)}{dx_i} D^\alpha\uvec \cdot D^\alpha\uvec \,d\xvec \bigg|
\leq c (1+ \eps  \|\uvec  \|_{s} ) 
\| D^\alpha\uvec \|_{0}^2 
\leq c (1+ \eps  \|\uvec  \|_{s} ) 
\|  \uvec \|_{s} ^2.
\end{gather}

The last term of (\ref{est1}) is left to estimate .
For this, we first notice that $D^\alpha  F_0$
is the sum of controlled coefficients 
multiplied by $D^\beta\uvec$ (possibly multiplied by $\eps^2$ or $\eps^4$ or \dots )
%or $\eps^2D^\beta\phhdet$ 
with 
$\beta = (\beta_1,\beta_2) \in \nit^2$ such that $\beta\leq\alpha$ (i.e.
$\beta_1\leq\alpha_1$ and $\beta_2\leq\alpha_2$). 
Hence  
\begin{gather}
\label{EstIneg5}
\bigg| \int D^\alpha  F_0 
\cdot D^\alpha\uvec \,d\xvec \bigg| 
\leq c (1+\sum_\beta \|D^\beta\uvec\|_0) \, \|D^\alpha\uvec\|_0
\leq c(1+\|\uvec \|_{s})\,\|\uvec \|_{s}.
\end{gather}
Secondly  
$[D^\alpha , A^1\fracp{}{x_1} ]\uvec$ is the sum of controlled coefficients 
multiplied by 
%\eps^{|\beta|} 
$D^\beta\uvec$ (possibly multiplied by $\eps$ or $\eps^2$ or \dots ) and themselves multiplied by
$D^\gamma\uvec$ with $\beta\leq \alpha$, $\gamma\leq\alpha$ and 
$\beta+\gamma\leq \alpha+(1,0)$ which implies $|\beta|+|\gamma|\leq |\alpha|+1$.
When $|\beta|\leq s-1$ and $|\gamma |\leq s-1$ since 
$|\beta|+|\gamma|+1\leq |\alpha|+2 < 2s$ we deduce that 
$D^\beta\uvec \cdot D^\gamma\uvec\in L^2(\rit^2)$  with
$\|D^\beta\uvec \cdot D^\gamma\uvec\|_0\leq c \|u\|_s^2 $ by a classical
calculus inequality that can be for instance found in the Appendix of 
Schochet \cite{schochet11:1986}.
When $|\alpha|=s$, $|\beta|= s$ and $|\gamma |=1$ we have 
$\sup_{\xvec\in\rit^2} |D^\gamma \uvec| \leq \|\uvec\|_s$. Then 
$D^\beta\uvec \cdot D^\gamma\uvec\in L^2(\rit^2)$.
When $|\alpha|=s$, $|\beta|= 1$ and $|\gamma |=s$  
$\sup_{\xvec\in\rit^2} |D^\beta \uvec| \leq \|\uvec\|_s$ and then 
$D^\beta\uvec \cdot D^\gamma\uvec\in L^2(\rit^2)$.
As the same can be done with 
$[D^\alpha , A^2\fracp{}{x_2}]\uvec$,
we deduce 
\begin{multline}
\label{EstIneg6}
\bigg| \int \Big(
- [D^\alpha , A^1\fracp{}{x_1} ]\uvec
- [D^\alpha , A^1\fracp{}{x_1} ]\uvec  \Big)
\cdot D^\alpha\uvec \,d\xvec \bigg| 
\leq c (1+\|\uvec \|_{s}+\|\uvec \|_{s}^2) \, \|D^\alpha\uvec\|_0
\\
\leq c(1+\|\uvec \|_{s}+\|\uvec \|_{s}^2)\,\|\uvec \|_{s}.
\end{multline}
Finally, $[D^\alpha , A^0\fracp{}{t} ]\uvec $ is a sum of controlled coefficients 
multiplied by $\eps D^\beta\uvec$ (possibly multiplied by $\eps$ or $\eps^2$ or \dots )
%$\eps^{|\beta|} D^\beta\uvec$ 
and themselves multiplied by
$D^\gamma\fracp{\uvec}{t}$ with $\gamma<\alpha$ and $(0,0)<\beta\leq\alpha$,
that implies $|\beta|>0$ and $|\beta|+|\gamma|\leq |\alpha|$.
When $|\beta|\leq s-1$ and $|\gamma |\leq s-2$, since 
$|\beta|+|\gamma|+1\leq |\alpha|+1 < s+(s-1)$, applying classical
calculus inequalities, we obtain
$D^\beta\uvec \cdot D^\gamma\fracp{\uvec}{t}\in L^2(\rit^2)$, with
$\|D^\beta\uvec \cdot D^\gamma\fracp{\uvec}{t}\|_0
\leq c \|u\|_s^{1/2}\|\fracp{\uvec}{t}\|_{s-2}^{1/2} \|u\|_{s-1}^{1/2}\|\fracp{\uvec}{t}\|_{s-1}^{1/2}$
$\leq c \|u\|_s \|\fracp{\uvec}{t}\|_{s-1}$.
When $|\alpha|=s$, $|\beta|= s$ and $|\gamma |=0$ we have 
$\sup_{\xvec\in\rit^2} |D^\gamma \fracp{\uvec}{t}| \leq \|\fracp{\uvec}{t}\|_{s-1}$.
Then $D^\beta\uvec \cdot D^\gamma\fracp{\uvec}{t}\in L^2(\rit^2)$.
When $|\alpha|=s$, $|\beta|= 1$ and $|\gamma |=s-1$, we get  
$\sup_{\xvec\in\rit^2} |D^\beta \uvec| \leq \|\uvec\|_s$ and 
$D^\beta\uvec \cdot D^\gamma\fracp{\uvec}{t}\in L^2(\rit^2)$.
Hence, we deduce 
\begin{gather}
\label{EstIneg9}
\bigg| \int \Big(
 [D^\alpha , A^0\fracp{}{t} ]\uvec\big)
\cdot D^\alpha\uvec \,d\xvec \bigg| 
\leq c(1+\eps \Big\| \fracp{\uvec}{t} \Big\|_{s-1}\|\uvec \|_{s})\,\|\uvec \|_{s}.
\end{gather}

Using inequalities (\ref{EstIneg1}) - (\ref{EstIneg9})
and summing (\ref{est1}) for $\alpha\leq s,$ we obtain
\begin{gather}
\label{est2}
\sum_{|\alpha|\leq s}
\bigg|
\frac{\ds d 
\Big(\int  A^0 D^\alpha\uvec \cdot D^\alpha\uvec \,d\xvec \Big)}{dt}
\bigg|
\leq g_1(\|\uvec \|_{s}, \eps\|\fracp{\uvec}{t} \|_{s-1} ),
\end{gather}
for a function $g_1$ not depending on $\eps$.

Derivating system (\ref{DalSystHyp}) with respect to $t$, we get
\begin{multline}
\label{ddtDalSystHyp}
A^0\fracp{(D^\alpha\fracp{\uvec}{t})}{t} 
+ A^1 \fracp{(D^\alpha\fracp{\uvec}{t})}{x_1} 
+ A^2 \fracp{(D^\alpha\fracp{\uvec}{t})}{x_2}
+ \frac{1}{\eps} S^1 \fracp{(D^\alpha\fracp{\uvec}{t})}{x_1} 
+ \frac{1}{\eps} S^2 \fracp{(D^\alpha\fracp{\uvec}{t})}{x_2}
+ \frac{1}{\eps}(D^\alpha\fracp{\uvec}{t})^\perp
\\
= -\frac{dA^0}{dt}D^\alpha\fracp{\uvec}{t}
-\frac{dA^1}{dt}D^\alpha\fracp{\uvec}{x_1}
-\frac{dA^2}{dt}D^\alpha\fracp{\uvec}{x_2}
+\frac{dF^\alpha}{dt}.
\end{multline}
Every previously established estimate remains valid. Moreover,
we can set
that the entries $A^i_{kl}$ of $A^i$ for $i=1,2$ satisfy 
\begin{gather}
\label{ddtEstIneg1}
\Big| \frac{d(A^i_{kl}(t,\tsep,\xvec,\eps \uvec))}{dt} \Big|
\leq c (\frac1\eps + \eps \big \| \uvec  \big \|_{s} 
+ \eps \big \| \fracp{\uvec}{t}  \big \|_{s-1} ),
\end{gather}
and using, when it is necessary, the same classical calculus procedure as above,
we can show that  $\frac{dF^\alpha}{dt}$ is in $L^2(\rit^2)$ with
\begin{gather}
\label{ddtEstIneg2}
\Big\| \frac{dF^\alpha}{dt}\Big\|_0
\leq c (\frac1\eps + 1
+ \eps \big \| \uvec  \big \|_{s} 
+ \eps \big \| \uvec  \big \|_{s}^2
+ \eps \big \| \fracp{\uvec}{t}  \big \|_{s-1} 
+ \eps \big \| \fracp{\uvec}{t}  \big \|_{s-1}^2 ) 
\big \| \uvec  \big \|_{s}.
\end{gather}
Then, multiplying equation (\ref{ddtDalSystHyp})
by $2 \eps D^\alpha(\eps\fracp{\uvec}{t})$ and following the previous method,
we obtain
\begin{gather}
\label{ddtest2}
\sum_{|\alpha|\leq s-1}
\bigg|
\frac{\ds d 
\Big(\int  A^0 D^\alpha(\eps\fracp{\uvec}{t} )
\cdot D^\alpha(\eps\fracp{\uvec}{t}) \,d\xvec \Big)}{dt}
\bigg|
\leq g_2(\|\uvec \|_{s}, \eps\|\fracp{\uvec}{t} \|_{s-1} ),
\end{gather}
for a function $g_2$ not depending on $\eps$.

As a conclusion, estimates (\ref{est2})-(\ref{ddtest2}) 
together with the fact that 
\begin{equation}
\Big(\sum_{|\alpha|\leq s}\int A^0 D^\alpha\uvec
\cdot D^\alpha\uvec \,d\xvec \Big)^{1/2}
\end{equation}
is a norm equivalent to $\|~\|_s$, allows us to set that 
\begin{equation}
\label{est3333}
\frac{\ds d\Big(\|\uvec\|_s +\eps \Big\| \fracp{\uvec}{t}\Big\|_{s-1} \Big)}{dt} 
\leq g(\|\uvec\|_s +\eps \Big\| \fracp{\uvec}{t}\Big\|_{s-1})
\end{equation}
for a function $g$ not depending on $\eps$ and then that the time interval
on which the classical solution of (\ref{SM1cs3intr}), 
(\ref{SM1cs2intr}) and (\ref{Cini})  exits does not depend on $\eps$.

~

Finally, estimate (\ref{estimsol}) is a direct consequence of  (\ref{est3333}).
This ends the proof of Theorem \ref{THM1}.

\section{Asymptotic behavior: proof of Theorem \ref{THM2}}
\label{AsBe}
In order to deduce the asymptotic behavior as $\eps$ goes to 0 of
$(\phhde ,\pvmr )$ we use the method, developed in 
Tartar \cite{tartar:1977},
Fr\'enod \cite{frenod:1994}
and Fr\'enod and Hamdache \cite{frenod/hamdache:1996} 
and used in Fr\'enod and Sonnendr\"ucker \cite{frenod/sonnendrucker:1997,frenod/sonnendrucker:1999}, 
which consists in setting  a weak formulation of system 
(\ref{SM1cs3intr}) - (\ref{SM1cs2intr}) or of its equivalent 
form (\ref{DalSystHyp}).
Passing then to the limit using the weak$-*$ convergence allows us to set a constraint
equation. This constraint equation imposes a form to $(\phhde ,\pvmr )$.
Using test functions satisfying the constraint equation
allows us finally to deduce system (\ref{SysEff}).

~

Estimate (\ref{estimsol}) yields the weak$-*$ 
convergence of $\uvec=(\phhdet,\pvmrt)$ to  
$\uuvec=(\phhde ,\pvmr )$ in $L^\infty([0,T],$ $H^s(\rit^2)^3)$, 
up to a subsequence, as $\eps$ goes to $0$.

Multiplying symmetric hyperbolic system (\ref{SystHyp}) by 
test functions $\Psi(t,\xvec)$ being
regular, $\rit^3-$ valued,  and with compact support 
in $[0,T)\times \rit^2$, and integrating yields
\begin{multline}
\label{ffftoc1.0}
-\int_0^T \int_{\rit^2} \uvec \cdot \bigg(
\fracp{A^0 \Psi}{t} 
+ \frac{1}{\eps} \fracp{A^0 }{\theta} \Psi
+ \fracp{ A^1 \Psi}{x_1} 
+ \fracp{ A^2 \Psi}{x_2}
+ \frac{1}{\eps} S^1 \fracp{\Psi}{x_1} 
+ \frac{1}{\eps} S^2 \fracp{\Psi}{x_2}
+ \frac{1}{\eps} \Psi^\perp
                             \bigg)  dt d\xvec
\\
= \int_0^T \int_{\rit^2} A^0 F \cdot \Psi dt d\xvec
+ \int_{\rit^2} \uvec_0 \cdot A^0 \Psi(0,0,\cdot) d\xvec.
\end{multline}
Multiplying (\ref{ffftoc1}) by $\eps$ and passing to the limit yields,
since $A^0$ converges to $I$ and $\partial A^0/\partial \theta$
to 0,
\begin{gather}
\label{ffftoc1}
-\int_0^T \int_{\rit^2} \int_0^1\uuvec \cdot \bigg(
S^1 \fracp{\Psi}{x_1} 
+ S^2 \fracp{\Psi}{x_2}
+ \Psi^\perp d\theta 
                             \bigg) dt d\xvec
=0,
\end{gather}
which is the weak formulation of
\begin{gather}
\label{eqcont1} 
S^1 \fracp{\uuvec}{x_1} 
+ S^2 \fracp{\uuvec}{x_2} 
+ \uuvec^\perp =0,
\end{gather}
or
\begin{gather}
\label{eqcont2}
\nabla\cdot\pvmr =0,
~~~~
\pvmr^\perp
+ \nabla\phhde =0.
\end{gather}
From this constraint equation, we deduce the form of
$(\phhde,\pvmr)$ given by (\ref{ContrCons}).

~

For any regular function $\varphi$ we define the 
test function $\Psi$ satisfying the constraint equation by  
\begin{equation}
\label{ContrFctTest}
\Psi_1(t,\xvec)= \varphi(t,\xvec) ~~~
\Psi_2(t,\xvec)= - \fracp{\varphi}{x_2} (t,\xvec)~~~
\Psi_3(t,\xvec)=\fracp{\varphi}{x_1} (t,\xvec) .
\end{equation}
Using this function in (\ref{ffftoc1.0}) cancels terms containing $1/\eps$
factors.
Since $A^0$ converges to $I$, 
$1/\eps\, \partial A^0/\partial\theta$ weak$-*$ converges 
to 0, $A^1$ weak$-*$ converges to $\int_0^1 \vmr_1 d\theta \, I$
and $A^2$ weak$-*$ converges to $\int_0^1 \vmr_2 d\theta \, I$ passing to the limit yields 
\begin{multline}
\label{ffEffEq1}
-\int_0^T \int_{\rit^2}\uuvec \cdot \bigg(
\fracp{\Psi}{t}
+ \fracp{ (\int_0^1 \vmr_1 d\theta) \Psi}{x_1} 
+ \fracp{ (\int_0^1 \vmr_2  d\theta) \Psi}{x_2}
                             \bigg)dt d\xvec 
\\
= \int_0^T \int_{\rit^2}\int_0^1 F d\theta \cdot \Psi dt d\xvec 
+ \int_{\rit^2} \uvec_0 \cdot \Psi(0,\cdot) d\xvec .
\end{multline}
or, using expressions of $\uuvec$ and $F$,
\begin{multline}
\label{ffEffEq2}
-\int_0^T \int_{\rit^2}  \phhde \Big(
\fracp{\varphi}{t}
+ \fracp{ (\int_0^1 \vmr_1 d\theta) \varphi}{x_1} 
+ \fracp{ (\int_0^1 \vmr_2  d\theta) \varphi}{x_2}\Big)
\\
+\fracp{\phhde}{x_2} \Big(
\fracp{\fracp{\varphi}{x_2}}{t}
+\fracp{ (\int_0^1 \vmr_1 d\theta) \fracp{\varphi}{x_2}}{x_1} 
+\fracp{ (\int_0^1 \vmr_2  d\theta) \fracp{\varphi}{x_2}}{x_2}\Big)
\\
+\fracp{\phhde}{x_1} \Big(
\fracp{\fracp{\varphi}{x_1}}{t}
+\fracp{ (\int_0^1 \vmr_1 d\theta) \fracp{\varphi}{x_1}}{x_1} 
+ \fracp{ (\int_0^1 \vmr_2  d\theta) \fracp{\varphi}{x_1}}{x_2}\Big)dt d\xvec 
\\
=\int_0^T \int_{\rit^2}
-\Big(\fracp{(\int_0^1 \hhde d\theta)}{x_1} (-\fracp{\phhde}{x_2}) 
+ \fracp{(\int_0^1 \hhde d\theta)}{x_2} (\fracp{\phhde}{x_1}) 
+ \big(\fracp{(\int_0^1 \vmr_1 d\theta)}{x_1}
+\fracp{(\int_0^1 \vmr_2 d\theta)}{x_2}\big) \phhde \Big)\varphi
\\
-\Big(\int_0^1 \vwn_1d\theta - 
(\fracp{\int_0^1 \vmr_1 d\theta}{x_1}(-\fracp{\phhde}{x_2}) +
\fracp{\int_0^1 \vmr_1 d\theta}{x_2}(\fracp{\phhde}{x_1})\Big)
\fracp{\varphi}{x_2}
\\
+\Big(\int_0^1 \vwn_2d\theta - 
(\fracp{\int_0^1 \vmr_2 d\theta}{x_1}(-\fracp{\phhde}{x_2}) +
\fracp{\int_0^1 \vmr_2 d\theta}{x_2}(\fracp{\phhde}{x_1})\Big)
\fracp{\varphi}{x_1}
dt d\xvec 
\\
+ \int_{\rit^2} 
\phhdet_0 \varphi 
- (\pvmrt_0)_1\fracp{\varphi}{x_2}
+ (\pvmrt_0)_2\fracp{\varphi}{x_1}d\xvec .
\end{multline} 
We have here a weak formulation of (\ref{SysEff}).
Since this equation is linear, it is easy to show that its solution 
is unique. From this, we can finally deduce that the whole sequence $\uvec$
weak$-*$ converges to $\uuvec$ as $\eps\rightarrow 0$,
ending the proof.

\section{Conclusion and perspectives}
\label{conclu}

In this paper, we set out equations modeling the long term evolution of the
perturbation $\phhdet$ of the ocean free surface elevation  and
of the perturbation $\pvmrt$ of the velocity field.
Because of the tide wave, those models contain and generate oscillations
with high frequency.
If numerical simulations of near coastal ocean waters are needed, directly using
those models could be very expensive because of the 
oscillations. 
Nevertheless, the result given in Theorem \ref{THM2} , which says 
that $(\phhdet,\pvmrt)$ weak$-*$ converges to $(\phhde ,\pvmr )$,
suggests a way to use those models for numerical simulations 
of near coastal ocean waters.
As an intuitive interpretation of it we could say that
\begin{gather}
\phhdet(t,\xvec) \text{ is close to } \phhde(t,\xvec)
\text{ ~  and  ~ }
\pvmrt(t,\xvec) \text{ is close to } \pvmr(t,\xvec).
\end{gather}
Hence, since equations (\ref{SysEff}) - (\ref{SysEffCI}) neither
contain nor generate oscillations with frequency $1/\eps$, 
we can solve it using a numerical method involving a time step which does not
need to be small compared with $\eps$.
Hence solving (\ref{SysEff})-(\ref{SysEffCI}) 
and then reconstructing $(\phhde ,\pvmr )$ via (\ref{ContrCons})
in place of  solving
(\ref{SM1cs3intr})-(\ref{SM1cs2intr})-(\ref{Cini}) directly
may give good results in a shorter computational time.

~

Among the tasks listed at the beginning of this paper, we realized 
significant steps in the direction of building numerical methods in our
previous paper \cite{AilFreMon} and in  coastal ocean modeling over
long time periods in the present one.

The next step, we shall provide in a forthcoming paper, will consist
in using the models set out here in order to compute the ocean fields
in a real near coastal ocean area and to couple this to the numerical
method proposed in \cite{AilFreMon}. This will allow us to make forecasts
in a real coastal ocean region.

\bibliographystyle{plain}
\bibliography{biblio}

\end{document}